\documentclass[12pt]{article}
\usepackage{amsmath,amssymb,amsfonts}
\usepackage{color}

\numberwithin{equation}{section}
\newtheorem{thm}{Theorem}

\def\cwedge{\bigcirc\kern-1.07em\wedge\ }
\def\owedge{\bigcirc\kern-0.82em\wedge}

\begin{document}

\title{Hermitian  structures on the product of Sasakian manifolds}

\author{Jung Chan Lee $^{1}$,
 Jeong Hyeong Park$^{1}$\thanks{Department of Mathematics, Sungkyunkwan University,
 Suwon 440-746, Korea, e-mail: parkj@skku.edu}, and
 Kouei Sekigawa$^{2}$\thanks{Department of Mathematics, Faculty of Science,
 Niigata University, Niigata, 950-2181, Japan,
 e-mail: sekigawa@math.sc.niigata-u.ac.jp}}

\date{}

\maketitle

\begin{abstract}
We investigate the curvature properties of a two-parameter family of
Hermitian structures on the product of two Sasakian manifolds, as
well as intermediate relations. We give a necessary and sufficient
condition for a Hermitian structure belonging to the family to be
Einstein and provide {{concrete}} examples.\\



\end{abstract}

\section{Introduction }\label{sec1}

In \cite{T}, Tsukada worked on the isospectral problem with respect
to the complex Laplacian for a two-parameter family of Hermitian
structures on the Calabi-Eckmann manifold $S^{2p+1} \times S^{2q+1}$
including the canonical one. In this paper, we define a
two-parameter family of almost Hermitian structures on the product
manifold $\overline{M} = M \times M'$ of a $(2p+1)$-dimensional
Sasakian manifold $M$ and a $(2q+1)$-dimensional Sasakian manifold
$M'$ similarly to the method used in \cite{T}, and show that any
almost Hermitian structure on $\overline{M}$ belonging to the two
parameter family is integrable; thereby generalizing the result of
the Calabi-Eckmann manifold by Tsukada (\cite{T}, {{Proposition
3.1}}). We shall further {{give}} a necessary and sufficient
condition for a Hermitian manifold in the family to be Einstein
(Theorem 1), {{which is the main result of the present paper. In the
last section, we shall provide concrete examples of Einstein
Hermitian structures on the Calabi-Eckmann manifolds $S^{2p+1}
\times S^{2q+1}$ for all $(p,q)$ $(p,q \geqq 1)$, by making use of
Theorem 1.}} {{However, we may also see that anyone of these
examples can not be weakly $*$-Einstein}}. {{Therefore}}, we may
note that the fact ``any compact Einstein Hermitian surface is
weakly $*$-Einstein" is not valid for higher dimensional cases in
general. {{The outline of these circumstances will be also given in
the last section}}.


\section{Preliminaries}\label{sec2}

In this section, we prepare some fundamental tools which we need in
our arguments. Let $\overline{M}=(\overline{M}, \bar{J}, \bar{g})$
be a $2n (\ge 4)$-dimensional almost Hermitian manifold with almost
Hermitian structure $(\bar{J},\bar{g})$. We denote by
$\overline{\nabla}$, $\bar{R}$, $\bar{\rho}$ and $\bar{\tau}$ the
Riemannian connection, the curvature tensor, the Ricci tensor and
the scalar curvature of $\overline{M}$, respectively. The curvature
tensor is defined by
\begin{equation*}
        \bar{R}(\bar{X},\bar{Y})\bar{Z}
        =[\overline{\nabla}_{\bar{X}},\overline{\nabla}_{\bar{Y}}]\bar{Z}-\overline{\nabla}_{[\bar{X},\bar{Y}]}\bar{Z},
    \end{equation*}
for $\bar{X}, \bar{Y}, \bar{Z}\in\mathfrak{X}(\overline{M})$, where
$\mathfrak{X}(\overline{M})$ denotes the Lie algebra of all smooth
vector fields on $\overline{M}$. The Ricci $*$-tensor $\bar{\rho}^*$
of $\overline{M}$ is defined by
\begin{equation}
\begin{gathered}
      \bar{\rho}^*(\bar{X},\bar{Y})=\text{tr }(\bar{Z}\mapsto \bar{R}(\bar{X},\bar{J}\bar{Z})\bar{J}\bar{Y}) \\
   =\frac{1}{2}\text{tr }(\bar{Z}\mapsto
   \bar{R}(\bar{X},{{\bar{J}\bar{Y})\bar{J}\bar{Z})}}),
\end{gathered}
\end{equation}
{{for $\bar{X}$, $\bar{Y}$, $\bar{Z}\in \mathfrak{X}(\overline{M})$.
It is easily checked
       that the equality $\bar{\rho}^* = \bar{\rho}$ holds on $\overline{M}$
       if $\overline{M}$ is K\"{a}hler.
}} We denote by $\bar{\tau}^{*}$ the $*$-scalar curvature of
$\overline{M}$, which is the trace of the Ricci $*$-operator $Q^*$
defined by
$\bar{g}(Q^*\bar{X},\bar{Y})=\bar{\rho}^*(\bar{X},\bar{Y})$. A
4-dimensional almost Hermitian manifold is also called an almost
Hermitian surface. We note that, for any almost Hermitian surface
$\overline{M}$, the Ricci and Ricci $*$-tensor are related by
\begin{equation}\label{eq:identity_for4}
\bar{\rho}^*(\bar{X}, \bar{Y}) + \bar{\rho}^*(\bar{Y}, \bar{X}) -
\{\bar{\rho}(\bar{X}, \bar{Y}) +\bar{\rho}(\bar{J}\bar{X},
\bar{J}\bar{Y}) \}= \frac{\bar{\tau}^*
-\bar{\tau}}{2}\bar{g}(\bar{X}, \bar{Y}),
\end{equation}
for $\bar{X}$, $\bar{Y}$, $\bar{Z} \in \mathfrak{X}(\overline{M})$
\cite{E,GH,K}. A $2n$-dimensional almost Hermitian manifold
$(\overline{M}, \bar{J},\bar{g})$ is called a {\it weakly
$*$-Einstein} manifold if the equality $\bar{\rho}^*=
\displaystyle{\frac{\bar{\tau}^*}{2n}\bar{g}}$ holds on
$\overline{M}$. Especially, if $*$-scalar curvature $\bar{\tau}^*$
of a weakly $*$-Einstein manifold $\overline{M}$ is constant,
 then $\overline{M}$ is said to be {\it $*$-Einstein.} It is
known that there exist weakly $*$-Einstein manifolds which are not
$*$-Einstein \cite{Ko,O}.  We denote by $\bar{N}$ the Nijenhuis
tensor of the almost complex structure $\bar{J}$ defined by
\begin{equation}\label{eq:2.4}
\bar{N}(\bar{X}, \bar{Y}) = [\bar{J}\bar{X}, \bar{J}\bar{Y}] -
[\bar{X}, \bar{Y}]- \bar{J}[\bar{J}\bar{X}, \bar{Y}]-
\bar{J}[\bar{X}, \bar{J}\bar{Y}]
\end{equation}
for $\bar{X}$, $\bar{Y}\in\mathfrak{X}(\overline{M})$. {{ {{It is
well-known that the almost complex structure $\bar{J}$ is integrable
if and only if the Nijenhuis tensor $\bar{N}$ vanishes identically
on $\overline{M}$ \cite{N}. An almost Hermitian manifold
$(\overline{M}, \bar{J}, \bar{g})$ with integrable almost complex
structure $\bar{J}$ is called a Hermitian manifold. It is well-known
that the condition $\bar{N}=0$ and the following condition
{{\eqref{eq:sasaki333}}} {{are}}
 equivalent:}}
\begin{equation}\label{eq:sasaki333}
\begin{aligned}
\bar{g}((\overline{\nabla}_{\bar{X}}\bar{J})\bar{Y},
\bar{Z})-\bar{g}((\overline{\nabla}_{\bar{J}\bar{X}}\bar{J})\bar{J}\bar{Y},
\bar{Z})=0
\end{aligned}
\end{equation}
for any vector field $\bar{X}, \bar{Y},
\bar{Z}\in\mathfrak{X}(\overline{M})$.

Next, we give a brief review of Sasakian manifolds. An almost
contact metric manifold $M=(M, \varphi, \xi, \eta, g)$ is called {{a
contact metric manifold if it satisfies}}
\begin{equation*}
d\eta(X,Y) = g(X, \varphi Y),
\end{equation*}
for any $X$, $Y$ $\in\mathfrak{X}(M)$. Further, a normal contact
metric manifold is called a Sasakian manifold. It is well-known that
a Sasakian manifold is characterized as an almost contact metric
manifold satisfying the condition
\begin{equation}
(\nabla_X \varphi)Y=g(X, Y)\xi-\eta(Y)X,
\end{equation}
for any $X$, $Y$ $\in\mathfrak{X}(M)$ {{(\cite{B}, Theorem 6.3)}}.
On a $(2n+1)$-dimensional Sasakian manifold $(M, \varphi, \xi, \eta,
g)$, we have the following identities :
\begin{equation}\label{eq:ce544444}
\begin{aligned}
 {{ \nabla_{X}\xi = - \varphi X, \hspace{5mm} (\nabla_{X}\eta)(Y) = - g(\varphi X, Y),}} \\
 {{ R(X,Y)\xi = \eta(Y)X - \eta(X)Y,}} \\
 {{\rho(\xi,X) = 2n \eta(X)}},
\end{aligned}
\end{equation}
for any $X$, $Y$ $\in\mathfrak{X}(M)$ \cite{B}. For further
investigation, we also prepare the following curvature identity on a
Sasakian manifold (\cite{B}, pp. 107):
\begin{equation}\label{eq:ce42}
\begin{aligned}
R(&X, Y, \varphi Z, W)-R(\varphi Z, X, Y, W)\\
&=-g(X, Y)g(\varphi Z, W)-2g(Z, \varphi Y)g(X, W)+g(Z, \varphi
X)g(Y, W),
\end{aligned}
\end{equation}
for $X$, $Y$, $Z$, $W\in\mathfrak{X}(M)$. From \eqref{eq:ce42}, we
get
\begin{equation}\label{eq:ce43}
\begin{aligned}
\sum_i^{2n+1}R(X, Y, \varphi e_i, e_i)-\sum_i^{2n+1}R(\varphi e_i,
X, Y, e_i)=3g(\varphi X, Y),
\end{aligned}
\end{equation}
for {{any}} orthonormal basis $\{e_1, \cdots , e_{2n+1}\}$ of
$T_xM$, $x\in M$.  Then,
 the left-hand side of \eqref{eq:ce43} implies
\begin{equation}\label{eq:ce44}
\begin{aligned}
\sum_i^{2n+1}R(X, Y, \varphi e_i, e_i)-\sum_i^{2n+1}R(\varphi e_i,
X, Y, e_i)=2\sum_i^{2n+1}R(X, Y, \varphi e_i, e_i)+\sum_i^{2n+1}R(Y,
\varphi e_i, X, e_i).
\end{aligned}
\end{equation}
Thus, from \eqref{eq:ce43} and \eqref{eq:ce44}, we have
\begin{equation}\label{eq:ce45}
\begin{aligned}
2\sum_i^{2n+1}R(X, Y, \varphi e_i, e_i)+\sum_i^{2n+1}R(Y, \varphi
e_i, X, e_i)=3g(\varphi X, Y).
\end{aligned}
\end{equation}
From \eqref{eq:ce45}, we  get also
\begin{equation}\label{eq:ce46}
\begin{aligned}
2\sum_i^{2n+1}R(Y, X, \varphi e_i, e_i)+\sum_i^{2n+1}R(X, \varphi
e_i, Y, e_i)=3g(\varphi Y, X).
\end{aligned}
\end{equation}
Thus, from \eqref{eq:ce45} and \eqref{eq:ce46}, we have
\begin{equation}\label{eq:ce47}
\begin{aligned}
4\sum_i^{2n+1}&R(X, Y, \varphi e_i, e_i)+\sum_i^{2n+1}R(Y, \varphi
e_i, X, e_i)+\sum_i^{2n+1}R(\varphi e_i, X, Y, e_i)=6g(\varphi X,
Y),
\end{aligned}
\end{equation}
and hence,
\begin{equation}\label{eq:ce48}
\begin{aligned}
\sum_i^{2n+1}R(X, Y, e_i, \varphi e_i)=-2g(\varphi X, Y).
\end{aligned}
\end{equation}
From \eqref{eq:ce48}, we have
\begin{equation}\label{eq:ce49}
\begin{aligned}
\sum_i^{2n+1}R(X, \varphi Y, e_i, \varphi e_i)=-2\big{(}g(X,
Y)-\eta(X)\eta(Y)\big{)},
\end{aligned}
\end{equation}
for any tangent vector $X$, $Y\in T_{x}M$, $x\in M$. The equality
\eqref{eq:ce49} is useful in the calculation of the formula
\eqref{eq:ce50} in the next section.
\section{{Curvature formulas and main result }}\label{sec3}
{{In this section, we define a two parameter family of almost
Hermitian structures on the product of Sasakian manifolds and show
the integrability. Further, we give a necessary and sufficient
condition for a Hermitian structure belonging to the family to be
Einstein one.}} Let $(M, \varphi, \xi, \eta, g)$ (resp. $(M',
\varphi', \xi', \eta', g')$) be a $(2p+1)$-dimensional Sasakian
manifold (resp. a $(2q+1)$-dimensional Sasakian manifold). We denote
by $\nabla$, $R$ and $\rho$ (resp. $\nabla'$, $R'$ and $\rho'$) the
Riemannian connection, the curvature tensor and the Ricci tensor on
$M$ (resp. $M'$). Let $\overline{M}=M\times M'$ be {{the product
manifold}} of $M$ and $M'$.
   Then we define a Riemannian metric $\bar{g}=\bar{g}_{a,b}${{ $(a, b
\in \mathbb{R})$ on $\overline{M}$}} by
 \begin{equation}\label{eq:sasaki2}
{{\bar{g}_{a,b}}}= g +a(\eta\otimes\eta' + \eta'\otimes\eta)+(a^2
+b^2 -1)\eta'\otimes\eta' +g'
\end{equation}
\cite{T}. {{Further, we define an almost complex structure $\bar{J}
= \bar{J}_{a,b}$ $(a,b \in \mathbb{R}, b\not= 0)$ }}on
$\overline{M}$ as follows.
\begin{equation}\label{eq:sasaki1}
\begin{aligned}
\bar{J}_{a,b}(X+X')= &\varphi(X)-\{\frac{a}{b}\eta(X) +
\frac{a^2+b^2}{b}\eta'(X')\}\xi\\
&+\varphi'(X')+\{\frac{1}{b}\eta(X)+\frac{a}{b}\eta'(X')\}\xi',
\end{aligned}
\end{equation}
for any tangent vector $X$ of $M$ and any tangent vector $X'$ of
$M'$ \cite{T}. It is easily checked that ${\bar{J}}^2= -I$ holds and
{{ $(\bar{J},\bar{g})$ is an almost Hermitian structure}} on
$\overline{M}$. 
 {{Since $\mathfrak{X}(M)$ and $\mathfrak{X}(M')$ are
regarded as the Lie subalgebra of $\mathfrak{X}(\overline{M})$, we
may note that (3.1) and (3.2) are rewritten respectively as
follows:}} 

\begin{equation}\label{eq:sasaki21}
\begin{aligned}
\bar{g}&(X, Y)=g(X, Y),\quad \bar{g}(X, Y')=a\eta(X)\eta'(Y'),\\
\bar{g}&(X', Y')=g'(X', Y')+(a^2+b^2-1)\eta'(X')\eta'(Y'),
\end{aligned}
\end{equation}
\begin{equation}\label{eq:sasaki211}
\begin{aligned}
\bar{J}(X)= &\varphi(X)-\frac{a}{b}\eta(X)\xi+\frac{1}{b}\eta(X)\xi',\\
\bar{J}(X')=
&\varphi'(X')-\frac{a^2+b^2}{b}\eta'(X')\xi+\frac{a}{b}\eta'(X')\xi',
\end{aligned}
\end{equation}
for $X$, $Y\in\mathfrak{X}(M)$, $X'$, $Y'\in\mathfrak{X}(M')$. Let
$\overline{\nabla}$, $\overline{R}$ and $\overline{\rho}$ be the
Riemannian connection, {{the}} curvature tensor and {{the}} Ricci
tensor of $\overline{M}$, respectively, {{and $X, Y, Z, W$ (resp.
$X', Y', Z', W'$) be any smooth vector field}} on $M$ (resp. $M'$).
{{Then, from (3.3) and (3.4), by making use of (2.5) and {{(2.6)}},
we have the following:
\begin{equation}\label{eq:sasaki3}
\begin{aligned}
\bar{g}(\overline{\nabla}_XY, Z)=&g(\nabla_XY, Z),\quad
\bar{g}(\overline{\nabla}_{X'}Y, Z)=-a\eta'(X')g(\varphi Y, Z),\\
\bar{g}(\overline{\nabla}_XY', Z)=&-a\eta'(Y')g(\varphi X, Z),\quad
\bar{g}(\overline{\nabla}_{X}Y, Z')=a\eta(\nabla_XY)\eta'(Z'),\\
\bar{g}(\overline{\nabla}_{X'}Y',
Z)=&a\eta'({\nabla'}_{X'}Y')\eta(Z),\quad
\bar{g}(\overline{\nabla}_{X'}Y, Z')=-a\eta(Y)g'(\varphi' X', Z'),\\
\bar{g}(\overline{\nabla}_XY', Z')=&-a\eta(X)g'(\varphi' Y',
Z'),\\
\bar{g}(\overline{\nabla}_{X'}Y', Z')=&g'({\nabla'}_{X'}Y', Z')
+(a^2+b^2-1)\big{(}\eta'({\nabla'}_{X'}Y')\eta'(Z')\\&-\eta'(X')g'(\varphi'
Y', Z')-\eta'(Y')g'(\varphi' X', Z')\big{)}.
\end{aligned}
\end{equation}
 By direct
calculation using {{ \eqref{eq:sasaki21},}} \eqref{eq:sasaki211} and
\eqref{eq:sasaki3}, we get the following:
\begin{equation}\label{eq:sasaki31111}
\begin{aligned}
\bar{g}((\overline{\nabla}_X\bar{J})Y, Z)=&\eta(Z)g(X,
Y)-\eta(Y)g(X,
Z),\\
\bar{g}((\overline{\nabla}_{X'}\bar{J})Y, Z)=&0,\\
\bar{g}((\overline{\nabla}_X\bar{J})Y', Z)=&b\eta'(Y')g(\varphi X,
Z)-a\eta'(Y')\big{(}g(X, Z)-\eta(X)\eta(Z)\big{)},\\
\bar{g}((\overline{\nabla}_{X'}\bar{J})Y',
Z)=&-a\eta(Z)\big{(}\eta'(X')\eta'(Y')-g'(X', Y')\big{)}+b\eta(Z)g'(\varphi' X', Y'),\\
\bar{g}((\overline{\nabla}_X\bar{J})Y', Z')=&0,\\
\bar{g}((\overline{\nabla}_{X'}\bar{J})Y',
Z')=&(a^2+b^2)\big{(}g'(X', Y')\eta'(Z')-g'(X', Z')\eta'(Y')\big{)}.
\end{aligned}
\end{equation}
Then, from {{(3.4)}} and (3.6), we can check that
\begin{equation}\label{eq:sasaki33}
\begin{aligned}
\bar{g}((\overline{\nabla}_{\bar{X}}\bar{J})\bar{Y},
\bar{Z})-\bar{g}((\overline{\nabla}_{\bar{J}\bar{X}}\bar{J})\bar{J}\bar{Y},
\bar{Z})=0
\end{aligned}
\end{equation}
holds for any $\bar{X}, \bar{Y},
\bar{Z}\in\mathfrak{X}(\overline{M})$.
Therefore, from (2.4), we see that the almost complex structure
$\bar{J}$ is integrable and hence, $(\overline{M},\bar{J}, \bar{g})$
is a Hermitian manifold, {{which generalizes the result of Tsukada
(\cite{T}, Proposition 3.1).}} However, from (3.6), we also see that
$\overline{M}$ is never K\"{a}hler. Now, we choose an orthonormal
basis $\{e_1, \cdots , e_{2p}, \xi\}$ (resp. $\{e'_1, \cdots ,
e'_{2q}, \xi'\}$) of the tangent space $T_xM$ (resp. $T_{x'}M'$) at
each point $x \in{M}$ (resp. $x' \in {M'}$). Then we may easily
check that $\{e_1, \cdots , e_{2p}, e'_1, \cdots , e'_{2q}, \xi,
 \frac{\xi'-a\xi}{b}\}$ is an orthonormal basis of
$T_{(x, x')}\overline{M}$,} which will be useful in the forthcoming
calculations of the present paper. {{From \eqref{eq:sasaki21} and
\eqref{eq:sasaki3}, by taking account of (2.5) and (2.6), we have
the formulas for the curvature tensor of $\overline{M}$: 
\begin{equation}\label{eq:ce16}
\begin{aligned}
\bar{g}(\bar{R}(X, Y)Z, W)=&g(R(X, Y)Z, W),\\
\bar{g}(\bar{R}(X, Y')Z, W)=&-a\eta'(Y')\big{(}g(X, Z)\eta(W)-g(X,
W)\eta(Z)\big{)},\\
\bar{g}(\bar{R}(X', Y')Z, W)=&2ag'(\varphi'X', Y')g(\varphi Z, W),\\
\bar{g}(\bar{R}(X, Y')Z, W')=&ag(\varphi X, Z)g'(\varphi'Y',
W')-a^2\eta'(Y')\eta'(W')\big{(}g(X,
Z)-\eta(X)\eta(Z)\big{)}\\
&-a^2\eta(X)\eta(Z)\big{(}g'(Y', W')-\eta'(Y')\eta'(W')\big{)},\\
\bar{g}(\bar{R}(X', Y')Z',
W)=&a(a^2+b^2)\eta(W)\big{(}\eta'(X)g'(Y', Z')-\eta'(Y')g'(X',
Z')\big{)},\\
\bar{g}(\bar{R}(X, Y)Z, W')=&a\eta(R(X, Y)Z)\eta'(W'),
\end{aligned}
\end{equation}
\begin{equation*}
\begin{aligned}
\bar{g}(\bar{R}(X', Y')Z', W')=g'(R'(X', Y')Z'&,
W')\\
+2(a^2+b^2-1)\big{\{}&\eta'(X')\eta'(W')g'(Y',
Z')-\eta'(Y')\eta'(W')g'(X', Z')\\
&-\eta'(X')\eta'(Z')g'(Y',
W')+\eta'(Y')\eta'(Z')g'(X', W')\big{\}}\\
-(a^2+b^2-1)^2\big{\{}&\eta'(X')\eta'(Z')g'(Y',
W')-\eta'(Y')\eta'(Z')g'(X', W')\\&+\eta'(Y')\eta'(W')g'(X',
Z')-\eta'(X')\eta'(W')g'(Y', Z')\big{\}}\\
+(a^2+b^2-1)\big{\{}&2g'(\varphi'X', Y')g'(\varphi'Z',
W')+g'(\varphi'X', Z')g'(\varphi'Y', W')\\
&-g'(\varphi'Y', Z')g'(\varphi'X', W')\big{\}}.
\end{aligned}
\end{equation*}
From \eqref{eq:ce16}, by direct {{calculation}}, we have
\begin{equation}\label{eq:ce18}
\begin{aligned}
\bar{\rho}(Y, Z)=&\rho(Y, Z)+2a^2q\eta(Y)\eta(Z),\\
\bar{\rho}(Y,
Z')=&2a\big{(}p+q(a^2+b^2)\big{)}\eta(Y)\eta'(Z'),\\
\bar{\rho}(Y', Z')=&\rho'(Y', Z')-2(a^2+b^2-1)g'(Y', Z')\\
&+2\big{(}pa^2+a^2+b^2-1+q(a^2+b^2-1)(a^2+b^2+1)\big{)}\eta'(Y')\eta'(Z'),
\end{aligned}
\end{equation}
for $Y$, $Z\in\mathfrak{X}(M)$, $Y'$, $Z'\in\mathfrak{X}(M')$. Thus,
from \eqref{eq:sasaki21} and \eqref{eq:ce18}, we see that
$(\overline{M}, \bar{g})$ is Einstein if and only if there is a
constant $\lambda$ satisfying the following conditions:
\begin{equation}\label{eq:ce19}
\begin{aligned}
\rho(Y, Z)+2a^2q\eta(Y)\eta(Z)=\lambda g(Y, Z),
\end{aligned}
\end{equation}
\begin{equation}\label{eq:ce20}
\begin{aligned}
2a\big{(}p+q(a^2+b^2)\big{)}\eta(Y)\eta'(Z')=a\lambda\eta(Y)\eta'(Z'),
\end{aligned}
\end{equation}
\begin{equation}\label{eq:ce21}
\begin{aligned}
&\rho'(Y', Z')-2(a^2+b^2-1)g'(Y', Z')\\
&+2\big{(}pa^2+a^2+b^2-1+q(a^2+b^2-1)(a^2+b^2+1)\big{)}\eta'(Y')\eta'(Z')\\
&=\lambda\big{(}g'(Y', Z')+(a^2+b^2-1)\eta'(Y')\eta'(Z')\big{)},
\end{aligned}
\end{equation}
for $Y$, $Z\in\mathfrak{X}(M)$, $Y'$, $Z'\in\mathfrak{X}(M')$. We
see that \eqref{eq:ce19} and \eqref{eq:ce21} may be rewritten as the
following, respectively.
\begin{equation}\label{eq:ce22}
\begin{aligned}
\rho(Y, Z)=\lambda g(Y, Z)-2a^2q\eta(Y)\eta(Z),
\end{aligned}
\end{equation}
\begin{equation}\label{eq:ce23}
\begin{aligned}
\rho'(Y', Z')=\big{(}\lambda&+2(a^2+b^2-1)\big{)}g'(Y', Z')\\
+\big{\{}&\lambda(a^2+b^2-1)-2\big{(}pa^2+a^2+b^2-1\\
&+q(a^2+b^2-1)(a^2+b^2+1)\big{)}\big{\}}\eta'(Y')\eta'(Z'),
\end{aligned}
\end{equation}
for $Y$, $Z\in\mathfrak{X}(M)$, $Y'$, $Z'\in\mathfrak{X}(M')$. Here,
by the assumption that $M$ and $M'$ are both Sasakian, we get
\begin{equation}\label{eq:ce24}
\begin{aligned}
\rho(Y, \xi)&=2p\eta(Y),\\
\rho'(Y', \xi')&=2q\eta'(Y'),
\end{aligned}
\end{equation}
for $Y\in\mathfrak{X}(M)$, $Y'\in\mathfrak{X}(M')$. Thus, from
\eqref{eq:ce22} and \eqref{eq:ce24}, we obtain
\begin{equation}\label{eq:ce25}
\begin{aligned}
2p\eta(Y)=(\lambda-2a^2q)\eta(Y).
\end{aligned}
\end{equation}
Similarly from \eqref{eq:ce23} and \eqref{eq:ce24}, we have also
\begin{equation}\label{eq:ce26}
\begin{aligned}
2q\eta'(Y')=\big{(}\lambda(a^2+b^2)-2pa^2-2q(a^2+b^2-1)(a^2+b^2+1)\big{)}\eta'(Y').
\end{aligned}
\end{equation}
Thus, from \eqref{eq:ce25} and \eqref{eq:ce26}, we obtain
\begin{equation}\label{eq:ce27}
\begin{aligned}
\lambda=2p+2a^2q,
\end{aligned}
\end{equation}
and
\begin{equation}\label{eq:ce28}
\begin{aligned}
(a^2+b^2)\lambda=2pa^2+2q(a^2+b^2)^2.
\end{aligned}
\end{equation}
Thus, from \eqref{eq:ce27} and \eqref{eq:ce28}, we get
\begin{equation}\label{eq:ce29}
\begin{aligned}
(p+a^2q)(a^2+b^2)=pa^2+q(a^2+b^2)^2,
\end{aligned}
\end{equation}
and hence
\begin{equation}\label{eq:ce30}
\begin{aligned}
pb^2+(a^2+b^2)a^2q=q(a^2+b^2)^2.
\end{aligned}
\end{equation}
From \eqref{eq:ce30}, we have further
\begin{equation}\label{eq:ce31}
\begin{aligned}
\big{(}p-(a^2+b^2)q\big{)}b^2=0.
\end{aligned}
\end{equation}
If $b\not=0$, then from \eqref{eq:ce31}, we have
\begin{equation}\label{eq:ce32}
\begin{aligned}
p=(a^2+b^2)q.
\end{aligned}
\end{equation}
Further, we suppose that $a\not=0$ (under $b\not=0$), then, from
\eqref{eq:ce20}, we have
\begin{equation}\label{eq:ce33}
\begin{aligned}
\lambda=4p.
\end{aligned}
\end{equation}
Thus, in this case, from \eqref{eq:ce27} and \eqref{eq:ce33}, we
have
\begin{equation}\label{eq:ce34}
\begin{aligned}
a^2q=p,
\end{aligned}
\end{equation}
and hence, taking account of \eqref{eq:ce32}, {{we have}}
\begin{equation}\label{eq:ce35}
\begin{aligned}
b^2=0.
\end{aligned}
\end{equation}
But, this is a contradiction. So, it must follow that $a=0$. Thus,
from \eqref{eq:sasaki21}, we have
\begin{equation}\label{eq:ce36}
\begin{aligned}
\bar{g}(Y, Z')=0,
\end{aligned}
\end{equation}
and
\begin{equation}\label{eq:ce37}
\begin{aligned}
\bar{g}(Y', Z')=g'(Y', Z')+(b^2-1)\eta'(Y')\eta'(Z'),
\end{aligned}
\end{equation}
for $Y\in\mathfrak{X}(M)$, $Y'$, $Z'\in\mathfrak{X}(M')$. Further,
from \eqref{eq:ce30}, we have also
\begin{equation}\label{eq:ce38}
\begin{aligned}
p=b^2q.
\end{aligned}
\end{equation}
From \eqref{eq:ce27}, we get
\begin{equation}\label{eq:ce39}
\begin{aligned}
\lambda=2p.
\end{aligned}
\end{equation}
From \eqref{eq:ce22} and \eqref{eq:ce23}, taking account of
\eqref{eq:ce38} and \eqref{eq:ce39}, we have
\begin{equation}\label{eq:ce40}
\begin{aligned}
\rho(Y, Z)=2pg(Y, Z),
\end{aligned}
\end{equation}
\begin{equation}\label{eq:ce41}
\begin{aligned}
\rho'(Y', Z')=2(p+b^2-1)g'(Y', Z')-2(b^2-1)(q+1)\eta'(Y')\eta'(Z'),
\end{aligned}
\end{equation}
for $Y$, $Z\in\mathfrak{X}(M)$, $Y'$, $Z'\in\mathfrak{X}(M')$.
Further, from \eqref{eq:ce18}, we have also
\begin{equation}\label{eq:ce411}
\begin{aligned}
\bar{\rho}(Y, Z')=0,
\end{aligned}
\end{equation}
for $Y\in\mathfrak{X}(M)$, $Z'\in\mathfrak{X}(M')$.\\
{{From \eqref{eq:ce544444},}} we can easily check that the Einstein
constant of any $(2p+1)$-dimensional Einstein Sasakian manifold is
equal to $2p$. Therefore, summing up the arguments above, we have
the following.
\begin{thm}\label{thm:2}
Let $M=(M, \varphi, \xi, \eta, g)$ and $M'=(M', \varphi', \xi',
\eta', g')$ be a $(2p+1)$-dimensional and a $(2q+1)$-dimensional
Sasakian manifold respectively, and let $\overline{M}=M\times M'$ be
the product manifold of $M$ and $M'$. {{Then}} $\overline{M}=
(\overline{M}, \bar{J},\bar{g})$ is a Hermitian manifold equipped
with the Hermitian structure $(\bar{g}, \bar{J})$ defined by
\eqref{eq:sasaki2} and \eqref{eq:sasaki1}. {{Furthermore}},
$\overline{M}=(\overline{M}, \bar{g}, \bar{J})$ is Einstein if and
only if $a=0$ and $M$ is an Einstein Sasakian manifold and $M'$ is
{{an}} $\eta$-Einstein Sasakian manifold with the Ricci tensor
$\rho'=2\displaystyle{(p+
\frac{p}{q}-1)g'-2(\frac{p}{q}-1)(q+1)\eta'\otimes\eta'}$.
\end{thm}


\noindent{\bf{Remark 1 }} We see that the $\eta$-Einstein Sasakian
manifold $M'$ {{in Theorem 1}} is Einstein if $p = q$. Then
{{$\overline{M}$}} is the Riemannian product of the
same dimensional Einstein Sasakian manifolds $M$ and $M'$. \\

Next, we shall calculate the Ricci $*$-tensor $\bar{\rho}^*$ of
$\overline{M}=(\overline{M}, \bar{g}, \bar{J})$. 
{{From \eqref{eq:sasaki211} and \eqref{eq:ce16}, by making use of
\eqref{eq:ce49}}}, we have
\begin{equation}\label{eq:ce50}
\begin{aligned}
\bar{\rho}^*(X, Y)=&(1-2aq)\big{(}g(X,
Y)-\eta(X)\eta(Y)\big{)},\\
\bar{\rho}^*(X,
Y')=&0,\\
\bar{\rho}^*(X',
Y')=&\big{(}1-2ap-(2q+1)(a^2+b^2-1)\big{)}\big{(}g'(X',
Y')-\eta'(X')\eta'(Y')\big{)},
\end{aligned}
\end{equation}
for $X$, $Y\in\mathfrak{X}(M)$, $X'$, $Y'\in\mathfrak{X}(M')$.\\

{{\noindent{\bf{{{Remark 2}}}} From \eqref{eq:sasaki21} and
\eqref{eq:ce50}, we see that $\overline{M}=(\overline{M}, \bar{J},
 \bar{g})$ is never weakly $*$-Einstein. }}

\section{Examples }\label{sec2}
Let $S^{2p+1}$ be a $(2p+1)$-dimensional unit sphere equipped with
the canonical Sasakian structure $(\varphi, \xi, \eta, g)$ of
constant sectional curvature $1$. Then, it can be seen that
$S^{2p+1}\times S^{2p+1}=(S^{2p+1}\times S^{2p+1}, \bar{J}_{0,1},
\bar{g}_{0,1})$ is {{an Einstein Hermitian }}manifold (See Remark
1). In this section, we provide further concrete examples of
{{Einstein Hermitian}} manifolds different from the above
{{trivial}} one based on the result of Theorem 1. A Sasakian
manifold with constant $\varphi$-holomorphic sectional curvature is
called a Sasakian space form.  First, we recall some fundamental
facts concerning a Sasakian space form. Let $M'=(M', \varphi', \xi',
\eta', g')$ be a $(2q+1) (\ge5)$-dimensional Sasakian space form
with constant $\varphi$-holomorphic sectional curvature $c$. Then,
it is known that the curvature tensor $R'$ is given by
\begin{equation}\label{eq:ce51}
\begin{aligned}
R'(X', Y')Z'=\frac{c+3}{4}\big{(}&g'(Y', Z')X'-g'(X', Z')Y'\big{)}\\
+\frac{c-1}{4}\big{(}&\eta'(X')\eta'(Z')Y'-\eta'(Y')\eta'(Z')X'\\
&+g'(X', Z')\eta'(Y')\xi'-g'(Y', Z')\eta'(X')\xi'\\
&+g'(\varphi'Y', Z')\varphi'X'-g'(\varphi'X', Z')\varphi'Y'\\
&-2g'(\varphi'X', Y')\varphi'Z'\big{)},
\end{aligned}
\end{equation}
for $X'$, $Y', Z'\in\mathfrak{X}(M')$ \cite{B,Ta}. Then, from
\eqref{eq:ce51}, we get easily
\begin{equation}\label{eq:ce52}
\begin{aligned}
\rho'(Y', Z')=&\frac{1}{2}\big{(}q(c+3)+c-1\big{)}g'(Y', Z')
-\frac{q+1}{2}(c-1)\eta'(Y')\eta'(Z'),
\end{aligned}
\end{equation}
for $Y'$, $Z'\in\mathfrak{X}(M')$.
Now, we shall introduce the following fact (\cite{B}, pp 114).
\begin{thm}\label{thm4}
Let $(M, \varphi, \xi, \eta, g)$ be a $(2q+1)(q>1)$-dimensional
Sasakian space form with constant $\varphi$-holomorphic sectional
curvature $c$ and apply the following D-homothetic deformation
\begin{equation*}\label{eq:ce5321}
\begin{aligned}
\eta'=\alpha\eta, \quad \xi'=\frac{1}{\alpha}\xi, \quad
\varphi'=\varphi,\quad g'=\alpha g+\alpha(\alpha-1)\eta\otimes\eta,
\end{aligned}
\end{equation*}
$\alpha$ being a positive constant, to the Sasakian structure
$(\varphi, \xi, \eta, g)$. Then, $(M, \varphi', \xi', \eta', g')$ is
a Sasakian space form with constant $\varphi$-holomorphic sectional
curvature $c'=\frac{c+3}{\alpha}-3$.
\end{thm}

{{Let $(S^{2q+1}, \varphi, \xi, \eta, g)$ be a $(2q+1)$-dimensional
unit sphere equipped with the canonical Sasakian structure
$(\varphi, \xi, \eta, g)$ of constant sectional curvature 1. Then,
by applying $D$-homothetic deformation with $\alpha=\frac{4}{c+3}
 (c>-3)$ to the canonical Sasakian structure $(\varphi, \xi, \eta,
g)$, we may obtain a Sasakian space form $(S^{2q+1}, \varphi', \xi',
\eta', g')$ of constant $\varphi$-holomorphic sectional curvature
$c$  by virtue of Theorem \ref{thm4}}}. Now, for any positive
integer $p, q (p\not=q)$, we set
\begin{equation}\label{eq:ce5421}
\begin{aligned}
c=\displaystyle{\frac{4p}{q}-3}.
\end{aligned}
\end{equation}
Then, from \eqref{eq:ce5421}, we may easily check that the both of
following equalities
\begin{equation}\label{eq:ce5321}
\begin{aligned}
\frac{1}{2}\big{(} q(c+3)+c-1\big{)}=2\displaystyle{(p+
\frac{p}{q}-1)},
\end{aligned}
\end{equation}
and
\begin{equation}\label{eq:ce54210}
\begin{aligned}
c-1=4\displaystyle{(\frac{p}{q}-1)}
\end{aligned}
\end{equation}
hold. Thus, from Theorem 1 and {{\eqref{eq:ce38},}} \eqref{eq:ce52},
\eqref{eq:ce5321}, \eqref{eq:ce54210}, we see that $(S^{2p+1}\times
S^{2q+1}, \bar{J}_{0,\sqrt{\frac{p}{q}}},
\bar{g}_{0,\sqrt{\frac{p}{q}}})$ $(p \ne q, q>1)$ {{provides a
non-trivial example of {{Einstein Hermitian manifold,}} where
{{$S^{2p+1}$ (resp. $S^{2q+1}$) is $(2p+1 )$-dimensional (resp.
$(2q+1)$-dimensional) unit sphere }} equipped with the canonical
Sasakian structures of constant sectional curvature $1$.
{{Therefore, taking account of {{the Remark}} 1, Remark 2 and the
examples provided in this section,}} {{we see that there exists a
$2n$-dimensional compact {{Einstein Hermitian}} manifold for any
integer $n (\geqq 3)$ which is not weakly $*$-Einstein.}} However,
the situation concerning this fact is quite different in the
4-dimensional case. In fact, by applying the Riemannian version of
the Goldberg-Sachs Theorem to a compact {{Einstein Hermitian}}
surface $M$ it follows that $M$ is {{a}} locally conformal K\"ahler
\cite{A,K}. We may easily check that the Ricci $*$-tensor $\rho^*$
is symmetric in a locally conformal K\"ahler surface since the Lee
form is closed {{(\cite{E}, (2.16)).}} Thus, taking account of the
identity {{\eqref{eq:identity_for4}}} on any almost Hermitian
      surface,  we see easily that $M$ is {{also}} weakly
$*$-Einstein. }}\\

{{\noindent{\bf{{{Remark 3}}}}  From (3.34),we see easily check that
the $*$-scalar curvature $\bar{\tau}^*$ of the Einstein Hermitian
manifold $(S^{2p+1}\times S^{2q+1}, \bar{J}_{0,\sqrt{\frac{p}{q}}},
\bar{g}_{0,\sqrt{\frac{p}{q}}})$ $(p \ne q, q>1)$ is given by
$\bar{\tau}^* = 4q(1 - p + q)$. Thus, taking account of {{the
Remark}} 1, we see that there exists a $2n$-dimensional compact
Einstein Hermitian manifold with constant $*$-scalar curvature for
any integer $n (\geqq 3)$ which is not K\"{a}hler, and this
completes the assertion of [6].}}

\section*{Acknowledgement}

This work was supported by the National Research Foundation of
Korea(NRF) grant funded by the Korea government(MEST)
(2011-0012987).

\end{document}